\documentclass{article}

\usepackage{amsfonts}
\usepackage{amssymb}
\usepackage{latexsym}

\newcommand{\pl}{\partial}

\newcommand{\be} {\begin{eqnarray}}
\newcommand{\ee} {\end{eqnarray}}
\newcommand{\ber} {\begin{eqnarray}}
\newcommand{\eer} {\end{eqnarray}}
\newcommand{\bep} {\begin{eqnarray*}}
\newcommand{\eep} {\end{eqnarray*}}

\newcommand{\dst} {\displaystyle}

\newcommand {\Hol}{\mathop{\rm Hol}\nolimits}

\renewcommand {\Re}{\mathop{\rm Re}\nolimits}

\newcommand {\GN}{\mathcal{GN}}
\newcommand {\NN}{\mathcal{N}}

\newcommand {\DD}{\Delta}

\newfont{\bbb}{msbm10}
\newfont{\bbbb}{msbm10 at 8pt}
\newfont{\bbbbb}{msbm10 at 11pt}

\def\Bbb#1{\hbox{\bbb #1}}

\newcommand{\R}{{\Bbb R}}

\newcommand{\B}{{\Bbb B}}

\newcommand{\C}{{\Bbb C}}

\newcommand {\HH}{{\cal H}}

\newcommand {\Ss}{{\cal S}}

\newcommand {\G}{{\cal G}}

\newtheorem{defin}{Definition}

\begin{document}

\title{A Julia--Carath\'eodory theorem for
hyperbolically monotone mappings in the Hilbert ball}

\author{Mark Elin\\ Simeon Reich\\ {\small and}\\ David Shoikhet}
\date{}

\maketitle
%\author{}

\begin{abstract}
We establish a Julia--Carath\'eodory theorem and a boundary
Schwarz--Wolff lemma for hyperbolically monotone mappings in the
open unit ball of a complex Hilbert space.
\end{abstract}
%\maketitle

\setcounter{equation}{0}

\indent Let $\B$ be the open unit ball of a complex Hilbert space
$\HH$ with inner product $\langle\cdot,\cdot\rangle$ and norm
$\|\cdot\|$, and let $\rho:\B\times\B\mapsto\R^+$ be the
hyperbolic metric on $\B$ (\cite{GK-RS-84}, p. 98), i.e.,
\be\label{rho}
\rho(x,y)=\tanh^{-1}\sqrt{1-\sigma(x,y)},
\ee
where
\be\label{sigma}
\sigma(x,y)=\frac{(1-\left\|x\right\|^2)(1-\left\|y\right\|^2)}
{|1-\langle x,y\rangle|^2},\quad x,y\in\B.
\ee

We denote by $\NN_\rho$ the class of all those self-mappings
$F:\B\mapsto\B$ which are nonexpansive with respect to $\rho$
($\rho$-nonexpansive), i.e.,
\begin{equation}        %equa1.2
\rho\biggl(F(x),F(y)\biggr)\leq\rho\left(x,y\right).
\end{equation}

Note that the class $\NN_\rho$ properly contains the class
$\Hol(\B)$ of all holomorphic self-mappings of $\B$ (\cite{FT-VE,
GK-RS-84}).

\begin{defin}%def1.1
A family $\Ss=\{F(t)\}_{t\geq 0}$ of self-mappings of $\B$ is said
to be a one-parameter continuous semigroup (flow) on $\B$ if
\be\label{semi_prop1}
F(t+s)=F(t)\circ F(s),\quad t,s\geq 0,
\ee
and
\begin{eqnarray}\label{semi_prop2}        %equa1.3
\lim_{t\rightarrow 0^{+}}F(t) =I,
\end{eqnarray}
where $I$ is the restriction of the identity mapping of $\HH$ to
$\B$ and the limit is taken pointwise with respect to the strong
topology of $\HH$.
\end{defin}

We denote $F(t)x$, the value of $F(t)$ at $x\in\B$, by $F_t(x),\
t\ge0$.

\begin{defin}    %def1.2
A flow $\Ss=\{F(t)\}_{t\ge0}$ on $\B$ is said to be generated if
for each $x\in\B$, there exists the strong limit
\begin{eqnarray}                %equa1.4
f(x):=\lim_{t\rightarrow 0^{+}}\,\frac 1t\biggl(x-F_t(x)\biggr).
\end{eqnarray}

In this case the mapping $f:\B\mapsto\HH$ is called the
(infinitesimal) generator of $\Ss$.

If $f$ generates a flow of $\rho$-nonexpansive self-mappings of
$\B$, then we will write $f\in\GN_\rho(\B)$.
\end{defin}

The following result is established in \cite{RS-SD1}:

\vspace{2mm}

$\blacklozenge$ {\sl A semigroup $\Ss$ of holomorphic
self-mappings of $\B$ is differentiable with respect to the
parameter $t\geq 0$ (hence, generated by a holomorphic mapping) if
and only if it is locally uniformly continuous on $\B$, i.e., the
limit in Definition 1 is uniform on each $\rho$-ball in $\B$.}

\vspace{2mm}

Moreover, in this case (see \cite{H-R-S} and \cite{RS2}):

\vspace{2mm}

$\blacklozenge$ {\sl The  generator $f$ is holomorphic on $\B$,
and bounded and uniformly continuous on each $\rho$-ball in $\B$.}

\vspace{2mm}

The set of all holomorphic semigroup generators is denoted by
$\G\Hol(\B)$.

\vspace{2mm}

The classical Julia--Carath\'{e}odory theorem and the boundary
Schwarz--Wolff lemma play a crucial role in geometric function
theory (see, for example, \cite{CC-BM} and \cite{SJH-93}). In
particular, they can be effectively used in the study of the
asymptotic behavior of discrete and continuous dynamical systems.
In this context these celebrated results may be stated as follows:

\vspace{2mm}

$\blacklozenge$ {\sl Let $F$ be a holomorphic self-mapping of the
open unit disk $\Delta$ in the complex plane $\C.$ If for a
boundary fixed point $\tau\in\partial\Delta$
$(\lim\limits_{r\to1^-}F(r\tau)=\tau)$ the angular derivative
$\angle F'(\tau):=\angle \lim\limits_{z\to\tau}F'(z) =k$ exists
finitely, then
\[
\frac{|F(z)-\tau|^2}{1-|F(z)|^2}\leq
k\,\frac{|z-\tau|^2}{1-|z|^2}\,.
\]
$\bullet$ If $k\leq 1$ this inequality means that each horocycle
internally tangent to the unit circle $\partial \Delta $ at
$\tau\in\partial\Delta$ is $F$-invariant.

\noindent$\bullet$ This is indeed the case when $F$ has no fixed
point in $\Delta$ and $\tau\in\partial\Delta $ is its so-called
Denjoy--Wolff point, that is, $\tau $ is an attractive fixed point
for all orbits of} $F.$

We use the symbol $\angle\lim\limits_{z\to\tau}$ to denote the
limit in each non-tangential approach region (see, for example,
\cite{PC} and \cite{SJH-93}).

\vspace{2mm}

Sometimes the above statements are grouped together under the name
{\it the Julia--Wolff--Carath\'eodory theorem}. Higher dimensional
analogs can be found, for instance, in \cite{Ru, GK-RS-84, DS-89,
CC-BM}.

For holomorphic mappings on the open unit disk $\Delta $ in the
complex plane $\C$ (i.e., for the one dimensional case when
$\B=\Delta$), an infinitesimal version of the
Julia--Wolff--Carath\'{e}odory theorem was given in \cite{ES1}.
Namely, it was shown there that

\vspace{2mm}

$\blacklozenge$ {\sl $f\in\G\Hol(\Delta)$ has no null point in
$\Delta$ if and only if for some $\tau\in\partial\Delta$, the
angular derivative
\begin{equation}\label{der}
\angle f'(\tau):=\angle \lim_{z\rightarrow \tau }f'(z) =\alpha
\end{equation}
exists (finitely) with $\Re \alpha \geq 0$.

Moreover,  $\alpha$ is, in fact, real and if
$\mathcal{S}=\{F(t)\}_{t\geq 0}$ is the semigroup generated by
$f$, then
\be\label{1star}
\frac{|F_t(z)-\tau|^2}{1-|F_t(z)|^2}\leq \exp
{(-t\alpha)}\,\frac{|z-\tau|^2}{1-|z|^2}\,,\quad z\in\Delta,\
t\ge0.
\ee
The point $\tau $ is unique and a (globally) attractive sink point
of $\mathcal{S}$.}

\vspace{2mm}

It is worth mentioning that the original Julia--Carath\'{e}odory
theorem deals not only with attractive boundary fixed points, but
also with repelling fixed points (see \cite{CC-BM, SR-DS-Enc}),
i.e., it deals with not necessarily fixed point free holomorphic
self-mappings. In this direction, a generalization of the above
theorem has recently been given by M.~D.~Contreras,
S.~D\'{i}az-Madrigal and Ch.~Pommerenke \cite{CM-DMS-CP}. They
proved the following one-dimensional assertion:

\vspace{2mm}

$\blacklozenge$ {\sl Let $f\in\G\Hol(\Delta).$ For a boundary
point $\tau \in
\partial \Delta $, the following claims are equivalent:

(i) the angular limit $\angle
f'(\tau):=\angle\lim\limits_{z\to\tau} \frac{f(z)}{z-\tau}=\beta$
exists finitely;

(ii) the angular limit $\frac{dF_t(\tau)}{dz}:=
\angle\lim\limits_{z\to\tau}\frac{F_t(z)-\tau}{z-\tau}=\exp(-t\beta)
.$}

\vspace{2mm}

Note that condition (ii) is equivalent to inequality (\ref{1star})
for some $\alpha \leq \beta .$

\vspace{2mm}

In this paper we will establish these assertions for a general
complex Hilbert (not necessarily finite-dimensional) space $\HH$.
Moreover, we will show that replacing the angular derivatives by
just radial derivatives, we are able to prove infinitesimal
versions of the Julia--Carath\'{e}odory theorem and the boundary
Schwarz--Wolff lemma for the much wider class of generators of
semigroups of $\rho$-nonexpansive self-mappings of the Hilbert
ball $\B$. For the case of the Denjoy--Wolff point, the asymptotic 
behavior of one-parameter semigroups of $\rho$-nonexpansive and holomorphic
mappings was also studied in \cite{E-R-S2}.

Indeed, the content of the classical Schwarz--Pick lemma is the
fact that each holomorphic self-mapping of the open unit disk in the
complex plane is nonexpansive with respect to the
Poincar\'e hyperbolic metric $\rho$ defined by~(\ref{rho}).
Therefore, one can try to use metric fixed point theory to derive
results regarding those mappings which are nonexpansive with
respect to $\rho$. At the same time, we must remember that if a
given mapping (or semigroup) is not holomorphic, then the notion
of derivative makes no sense.

It turns out, however, that although the study of the asymptotic
behavior of a fixed point free semigroup consisting of
$\rho$-nonexpansive mappings is, in general, much more complicated,
one can define the real part of the radial derivative of its generator at
a boundary fixed point and use it to find invariant ellipsoids,
internally tangent at this point to the unit sphere, in the spirit
of the Julia--Carath\'eodory theorem and the boundary
Schwarz--Wolff lemma. To understand this phenomenon, we first
consider the following example.

\vspace{2mm}

\textbf{Example 1.  }Consider the continuous function
$f:\,\overline{\Delta}\mapsto\C $ defined by
\[
f(z)=\frac{z+\bar{z}}{2}+\chi \,\frac{z-\bar{z}}{2}-1,
\]
where $\chi$ is a real parameter. Elementary calculations show
that for all $\chi\geq \frac12$ the following boundary flow
invariance condition holds:
\[
\Re f(z)\overline{z}\geq 0,\quad z\in\partial\Delta.
\]
Therefore it follows from Martin's theorem \cite{MRH-73} (see also
\cite{RS-76}) that $f$ generates a semigroup of continuous
self-mappings of $\overline{\Delta}.$ Indeed, solving the Cauchy
problem
\be\label{2a}
\left\{
\begin{array}{l}
{\displaystyle\frac{\pl F_t(z)}{\pl t}}+f(F_t(z))=0,\vspace{3mm}\\
F_0(z)=z,\quad z\in\Delta,
\end{array}
\right.
\ee
we find
\[
F_{t}(z)=1-e^{-t}+e^{-t}\frac{z+\bar{z}}{2}+e^{-\chi t}\frac{z-\bar{z}}{2}%
\,,\quad z\in \overline{\Delta }.
\]

It is clear that $F_{t}(1)=1$ for all $t\geq 0$ and
$\lim\limits_{t\rightarrow \infty }F_{t}(z)=1$ for all $z\in
\overline{\Delta }.$

However, no more information on invariant sets of this semigroup
can be obtained in this way.

If $\chi =1,$ then for each $t\geq 0$, the function $F_{t}(z),$ as
well as $f$, are holomorphic. Hence one can apply the result in
\cite{ES1} (see (\ref{der})--(\ref{1star}) above) to derive
\be\label{N}
d_{1}(F_{t}(z)):=\frac{|1-F_{t}(z)|^{2}}{1-|z|^{2}}\leq 
\frac{e^{-t}|1-z|^{2}%
}{1-|z|^{2}}\,,
\ee
because the angular derivative of $f$ at the boundary fixed point
$\tau =1$ exists and equals $1$.

But for $\chi >1$ our generator $f$, as well as its generated
semigroup, are again not holomorphic.

However, fortunately, in this situation ($\chi >1$) one can show
that the semigroup $\Ss=\left\{ F_{t}\right\}$ consists of
$\rho$-nonexpansive self-mappings of the open unit disk $\Delta$
with the Poincar\'e hyperbolic metric $\rho$ defined on it. Indeed,
\[
F_t(z)=\left(1-e^{-t}\right)+z\left(\frac{e^{-t}}2+\frac{e^{-\chi
t}}2 \right)+\bar{z}\left(\frac{e^{-t}}2-\frac{e^{-\chi t}}2
\right)
\]
is a convex combination of holomorphic and anti-holomorphic
(hence, $\rho$-nonexpansive) mappings. So, by Theorem~6.5 on
page~75 of \cite{GK-RS-84}, $F_t$ also is $\rho$-nonexpansive.

In addition, we see that although $f$ is not differentiable in the
complex sense in $\Delta,$ its radial derivative at the boundary
null point $\tau =1$ does exist:
\[
\alpha :=\lim_{r\rightarrow 1^{-}}\frac{f(r)}{r-1}=1.
\]

So, the following question arises:

$\bullet$\emph{\ Is this fact sufficient to ensure the same
invariance condition as (\ref{N}) for a semigroup of
$\rho$-nonexpansive self-mappings of the open unit disk $\Delta$?}

In this paper we give an affirmative answer to this question in a
much more general situation.

To formulate our results, we need the following notions and
notations.

\vspace{2mm}

For a fixed $\tau\in\pl\B$, the boundary of $\B$, and an arbitrary
$x\in\B$, we define a non-Euclidean ``distance" between $x$ and
$\tau $ by the formula
\begin{equation}
d_{\tau }(x)=\frac{\left\vert 1-\left\langle x,\tau \right\rangle
\right\vert ^{2}}{1-\left\Vert x\right\Vert ^{2}}\,.
\end{equation}

The sets
\be\label{ellips}
E(\tau,s)=\biggl\{x\in\B:\ d_\tau(x)=\frac{\left\vert
1-\left\langle x,\tau \right\rangle \right\vert ^{2}}{1-\left\Vert
x\right\Vert ^{2}}<s\biggr\},\quad s>0,
\ee
are ellipsoids internally tangent to the unit sphere $\partial\B$
at $\tau$.

As in \cite{AERS}, it can be shown  that the support functional
$x^{*}$ of the smooth convex set $E(\tau,s)$ at $x\in\partial
E(\tau,s),\ x\neq\tau$, normalized by the condition
\[
\lim_{x\rightarrow \tau }\left\langle x-\tau,x^{*}\right\rangle
=1,
\]
can be expressed by
\begin{equation}\label{funct}    %equ2.2
x^*=\frac 1{1-\|x\|^2}\,x-\frac 1{1-\langle \tau ,x\rangle
}\,\tau.
\end{equation}

\vspace{2mm}

\noindent\textbf{Theorem.} {\it Let $\mathcal{S}$ be a semigroup
of $\rho $-nonexpansive self-mappings of the Hilbert ball $\B$,
generated by $f:\,\B\mapsto\HH$. Suppose that $f$ is uniformly
continuous on each $\rho$-ball in $\B$, and $\tau\in\partial\B$ is
a null point of $f$ in the sense that
$\displaystyle\lim\limits_{r\rightarrow 1^{-}}f(r\tau )=0$. The
following assertions are equivalent:

\vspace{1mm}

\ber\label{limsup}
(I) \hspace{33mm}
\limsup_{r\to1^-}\Re\frac{\left<f(r\tau),\tau\right>}{r-1}>-\infty;
\hspace{33mm}\eer

\noindent(II) \quad
$\alpha:=\dst\lim\limits_{r\to1^-}\Re\,\frac{\langle
f(r\tau),\tau\rangle}{r-1}$ exists finitely;

\noindent(III) \quad $\beta:=\inf2\Re\langle
f(x),x^*\rangle>-\infty$;

\noindent(IV) there exists a real number $\gamma$ such that
\ber\label{phi}
d_\tau(F_t(x))\le\exp(-t\gamma)\cdot d_\tau(x),\quad x\in\B.
\eer

Moreover,

(a) $\alpha=\beta$ and the maximal $\gamma$ for which (IV) holds
is exactly the same $\beta$;

(b) if $f$ is holomorphic and one (hence all) of conditions
(I)--(IV) holds, then
\[
\lim\limits_{r\rightarrow 1^{-}}\frac{\langle f(r\tau ),\tau \rangle }{%
r-1}
\]%
exists and is actually a real number.}

\vspace{2mm}

Combining this result with Theorem 8.3 in \cite{RS-SD2005}, we arrive
at the following analog of the boundary Schwarz--Wolff lemma.

\vspace{2mm}

\noindent\textbf{Corollary.} {\it Let $\Ss=\{F_{t}\}_{t\geq
0}\subset\NN_\rho$ be a semigroup of $\rho$-nonexpansive
self-mappings of $\B$ generated by $f\in\GN_\rho(\B)$. Assume that
$f$ is null point free (that is, $\Ss$ has no stationary point in
$\B$). Then there is a unique point $\tau\in\partial\B$ such that
\[
d_{\tau }(F_{t}(x))\leq \exp (-t\alpha )\cdot d_{\tau }(x),\quad
x\in\B,
\]
for some $\alpha\ge0$, and there is a continuous curve $\{x(r):\
0\le r<1 \}\subset\B$ ending at $\tau\in\pl\B$ for which
\[
\lim_{r\rightarrow 1^-}f(x(r))=0.
\]
Conversely, if for some point $\tau\in\partial\B$ the radial limit
$\lim\limits_{r\to1^-}f(r\tau)=0$ and the radial limit $\alpha
:=\displaystyle\lim\limits_{r\rightarrow 1^{-}}\Re \frac{\langle
f(r\tau ),\tau \rangle }{r-1}$ exists and is positive, then for
each $t>0$, the mapping $F_{t}$ has no fixed point in~$\B.$}

\vspace{2mm}

To prove our theorem we need the following additional concepts and
facts.

\vspace{2mm}

A mapping $f:\B\mapsto\HH$ is said to be hyperbolically monotone
or $\rho$-monotone (\cite{RS2, K-R}) if for each pair $x,y\in \B$,
\begin{equation}
\rho \biggl(x+rf(x),y+rf(y)\biggr)\geq \rho (x,y)
\end{equation}%
for all $r\geq 0$ such that the points $x+rf(x)$ and $y+rf(y)$
belong to $\B$.

The crucial point in our approach is the following
characterization of $\rho$-monotonicity \cite[Theorem 2.1]{RS2}:

\vspace{2mm}

$\blacklozenge$ {\sl A mapping $f:\B\mapsto\HH$ is $\rho$-monotone
if and only if
\begin{equation}\label{rho_mon}
\mathop{\rm Re}\nolimits\left[ \frac{\langle f(x),x\rangle
}{1-\Vert x\Vert
^{2}}+\frac{\langle y,f(y)\rangle }{1-\Vert y\Vert ^{2})}\right] \geq %
\mathop{\rm Re}\nolimits\left[ \frac{\langle f(x),y\rangle
+\langle x,f(y)\rangle }{1-\langle x,y\rangle }\right]
\end{equation}
for all points $x,y\in\B$.}

Moreover, {\sl if $f\in\GN_\rho(\B)$ is uniformly continuous on
each $\rho$-ball in $\B$, then it is $\rho$-monotone.}

So, {\it each holomorphic generator on $\B$ is $\rho$-monotone.}

\vspace{2mm}

\noindent{\bf Proof of Theorem. } The implication
(II)$\Rightarrow$(I) is trivial. To prove other implications we
denote by $\psi$ the following real-valued function:
\[
\psi(t,x):=d_\tau(F_t(x)).
\]
By direct calculation we get
\be\label{2}
\left. \frac{\pl\psi(t,x)}{\pl
t}\right|_{t=0}=-2\psi(0,x)\Re\langle f(x),x^*\rangle,
\ee
where $x^*$ is defined by (\ref{funct}). Hence, by the semigroup
property,
\be\label{1}
\frac{\pl\psi(t,x)}{\pl t}=-2\psi(t,x)\Re\langle
f(F_t(x)),(F_t(x))^*\rangle
\ee
for all $t\ge0$.

\noindent{\bf Step 1.} (I)$\Rightarrow$(III). If (I) holds, then
there exists a sequence $\{r_n\}$ such that $r_n\nearrow1$ and
\[
\alpha_1:=\limsup_{r\to1^-}\Re\frac{\left<f(r\tau),\tau\right>}{r-1}
=\lim_{n\to\infty}\Re\frac{\langle
f(r_n\tau),\tau\rangle}{r_n-1}\,.
\]

Setting $y=r_n\tau$ in (\ref{rho_mon}), we get
\bep
\Re \frac{\langle f(x),x\rangle}{1-\|x\|^2}+ \frac{\langle
f(r_n\tau),r_n\tau\rangle}{1-r_n^2} \ge \Re \frac{\langle
f(x),r_n\tau\rangle+\langle x,f(r_n\tau)\rangle}{1-\langle
x,r_n\tau\rangle}\,,
\eep
or, equivalently,
\ber\label{wrtau}
\dst\Re\left<f(x),\frac{x}{1-\|x\|^2}-\frac{r_n\tau}{1-r_n\left<\tau,x\right>}\right>
\ge\nonumber\\[2mm]
\Re\left<f(r_n\tau),\frac{x}{1-r_n\left<x,\tau\right>}-\frac{r_n\tau}{1-r_n^2}\right>.
\eer
Letting now $n$ tend to $\infty$, we see that
\bep
\dst\Re\left<f(x),x^*\right> \ge \frac{\alpha_1}2\,,
\eep
i.e., (III) holds and $\beta$ is not less than $\alpha_1$. By the
way, this implies that $\alpha_1$ is finite.

\noindent{\bf Step 2.} (III)$\Leftrightarrow$(IV). Let (IV) hold
for some real number $\gamma$. Differentiating (\ref{phi}) at
$t=0$, we obtain
\[
2\Re\langle f(x),x^*\rangle\ge\gamma,
\]
i.e., (III) holds, and $\beta\ge\gamma$.

Let now (III) hold. Then by (\ref{1}),
\[
\frac{\pl\psi(t,x)}{\pl t}\le-\psi(t,x)\beta.
\]
Integrating this inequality with respect to $t$, we see that (IV)
holds with $\gamma=\beta$.

\noindent{\bf Step 3.} (III)$\Rightarrow$(I). Let
\[
2\Re \left<
f(x),\frac{x}{1-\|x\|^2}-\frac\tau{1-\langle\tau,x\rangle}\right>\ge\beta.
\]
Substituting $x=r\tau,\ 0<r<1,$ we see that
\[
2\Re \left<
f(r\tau),\frac{r\tau}{1-r^2}-\frac\tau{1-r}\right>\ge\beta,
\]
or
\[
\Re \frac{\langle
f(r\tau),\tau\rangle}{r-1}\cdot\frac2{r+1}\ge\beta.
\]
Therefore
\[
\liminf_{r\to1^-}\Re\frac{\langle
f(r\tau),\tau\rangle}{r-1}\ge\beta.
\]
Hence (I) holds and
$\limsup\limits_{r\to1^-}\Re\frac{\left<f(r\tau),\tau\right>}{r-1}\ge\beta$.

\noindent{\bf Step 4.} Just by comparing Step 1 and Step 3, we
conclude that (III) implies (II), i.e., if (III) holds, then
$\alpha=\lim\limits_{r\to1^-}\Re\frac{\left<f(r\tau),\tau\right>}{r-1}$
exists and is equal to $\beta$.

\noindent{\bf Step 5.} To end the proof, we have to prove (b).

Suppose that $f$ is holomorphic. We introduce a holomorphic
function $g\in\Hol(\DD,\C)$ as follows:
\ber\label{f1}
g(\lambda):=\left<f(\lambda\tau),\tau
\right>-\frac\beta2(\lambda^2-1).
\eer

It follows from (III) that
\[
\Re\left<f(\lambda\tau),(\lambda\tau)^*\right>\ge\frac\beta2\,.
\]
On the other hand,
\bep
&&\left<f(\lambda\tau),(\lambda\tau)^*\right>=
\left<f(\lambda\tau),\frac{\lambda\tau}{1-|\lambda|^2}-\frac{\tau}{1-\bar\lambda}\right>\\[2mm]
&&=\left<f(\lambda\tau),\tau\right>\cdot\Bigl(\frac{\bar\lambda}{1-|\lambda|^2}-\frac1{1-\lambda}\Bigr)
\\[2mm]
&&=\left(g(\lambda)+\frac\beta2(\lambda^2-1)\right)\cdot
\frac{\bar\lambda-1}{(1-|\lambda|^2)(1-\lambda)} \\[2mm]
&&=
g(\lambda)\cdot\frac{|\lambda-1|^2}{-(1-|\lambda|^2)(1-\lambda)^2}+\frac\beta2
\frac{(1-\bar\lambda)(1+\lambda)}{1-|\lambda|^2}\,.
\eep

Therefore $\ \dst\Re\frac{g(\lambda)}{-(1-\lambda)^2}\ge0$. By the
Riesz--Herglotz representation formula,
\[
\frac{g(\lambda)}{-(1-\lambda)^2}=\oint_{|\zeta|=1}
\frac{1+\lambda\bar\zeta}{1-\lambda\bar\zeta}\,d\sigma(\zeta),
\]
where $d\sigma(\zeta)$ is a positive measure on the unit circle.
Decomposing $d\sigma$ with respect to Dirac's $\delta$-function at
the point $\zeta=1$,
$d\sigma(\zeta)=a\delta(\zeta)+d\sigma_1(\zeta),\ a\ge0$, we
calculate:
\ber
\lim_{r\to1^-}\frac{g(r)}{r-1}=\lim_{r\to1^-}(1-r)
\oint_{|\zeta|=1}\frac{1+r\bar\zeta}{1-r\bar\zeta}\,(a\delta(\zeta)+d\sigma_1(\zeta))=
2a\ge0.
\eer

This fact, in turn, implies that
\ber
&& \lim_{r\to1^-}\frac{\left<f(r\tau),\tau\right>}{r-1}=
\lim_{r\to1^-}\left(\frac{g(r)}{r-1}+\frac\beta2(r+1)\right)=2a+\beta
\eer
exists and is real.

Since by Step 4, $\Re(2a+\beta)=\beta$, we have that $a=0$, and
hence the limit
\[
\lim_{r\to1^-}\frac{\left<f(r\tau),\tau\right>}{r-1}=\beta
\]
is real. This completes the proof of our theorem.

\vspace{2mm}

\noindent{\bf Example 2.} Let $\HH_1$ be a Hilbert space. 
Consider the semigroup of holomorphic self-mappings of the unit
ball $\B$ of the space $\HH:=\C\times\HH_1$ defined by the
formula
\bep
F_t(z_1,z_2)=\left(\frac{z_1}{z_1+e^t(1-z_1)}\,,\,\frac{e^{t/2}z_2}{z_1+e^t(1-z_1)}
\right),
\eep
where $z_1\in\C,\ z_2\in\HH_1,\ |z_1|^2+\|z_2\|^2<1$. This
semigroup has a boundary fixed point $\tau=(1,0)$ which is not its
Denjoy--Wolff point. Consider the semigroup generator
$f\in\Hol(\B,\HH)$:
\bep
f(z_1,z_2)=-\left.\frac{\pl F_t(z_1,z_2)}{\pl
t}\right|_{t=0}=\left(z_1(1-z_1),\frac{(1-2z_1)z_2}2\right).
\eep
Now we just calculate
\[
\alpha=\lim_{r\to1^-}\frac{\langle f(r\tau),\tau\rangle}{r-1}=-1.
\]
Hence,
\[
d_\tau\left(F_t(z)\right)\le e^t d_\tau(z),\quad z\in\B,\ t\ge0.
\]

\vspace{2mm}

\noindent{\bf Example 3.} Define another semigroup $\left\{F_t
\right\}_{t\ge0}=\left\{\left(F_t\right)_1,
\left(F_t\right)_2\right\}_{t\ge0}$ on the unit ball $\B$ of the
same space $\HH$ by the formulae:
\bep
&&\left(F_t\right)_1(z_1,z_2)=
\frac{(1+z_1^2)e^{2t}-(1-z_1)\sqrt{2(1+z_1^2)e^{2t}-(1-z_1)^2}}{(1+z_1^2)e^{2t}-(1-z_1)^2}\,,\\
&&\left(F_t\right)_2(z_1,z_2)=z_2\sqrt{\frac{\pl\left(F_t\right)_1(z_1,z_2)}{\pl
z_1}}\,,
\eep
where $z_1\in\C,\ z_2\in\HH_1,\ |z_1|^2+\|z_2\|^2<1$. Consider its
generator $f\in\Hol(\B,\HH)$:
\bep
f(z_1,z_2)=-\left.\frac{\pl F_t(z_1,z_2)}{\pl
t}\right|_{t=0}=\left(-\frac{(1-z_1)(1+z_1^2)}{1+z_1}\,,\frac{z_2(1-z_1+z_1^2+z_1^3)}{(1+z_1)^2}\right).
\eep
It is clear that $f$ has three boundary null points:
$\tau_1=(1,0),\ \tau_2=(i,0)$ and $\tau_3=(-i,0)$. For each one of
them we just calculate
\bep
&&\alpha_1=\lim_{r\to1^-}\frac{\langle
f(r\tau_1),\tau_1\rangle}{r-1}=\lim_{r\to1^-}\frac{-(1-r)(1+r^2)}{(1+r)(r-1)}=1,\\
&&\alpha_2=\lim_{r\to1^-}\frac{\langle
f(r\tau_2),\tau_2\rangle}{r-1}=\lim_{r\to1^-}\frac{-(1-ri)(1-r^2)(-i)}{(1+ri)(r-1)}=-2,\\
&&\alpha_3=\lim_{r\to1^-}\frac{\langle
f(r\tau_3),\tau_3\rangle}{r-1}=\lim_{r\to1^-}\frac{-(1+ri)(1-r^2)i}{(1-ri)(r-1)}=-2.
\eep
Hence, the following three inequalities hold simultaneously for
all $z\in\B$ and $t\ge0$:
\bep
d_{\tau_1}\left(F_t(z)\right)\le e^{-t} d_{\tau_1}(z),\\
d_{\tau_2}\left(F_t(z)\right)\le e^{2t} d_{\tau_2}(z),\\
d_{\tau_3}\left(F_t(z)\right)\le e^{2t} d_{\tau_3}(z).
\eep
These inequalities mean that for each point $z\in\B$ and for each $t\ge0$, 
the
image $F_t(z)$ belongs to the intersection of the ellipsoids:
\[
E(\tau_1,e^{-t} d_{\tau_1}(z))\bigcap E(\tau_2,e^{2t}
d_{\tau_2}(z))\bigcap E(\tau_3,e^{2t} d_{\tau_3}(z)).
\]

\vspace{2mm}

\noindent{\bf Example 4. }Consider the one-parameter continuous
semigroup $\Ss=\{F_t\}_{t\ge0}$, consisting of holomorphic
self-mappings of the open unit disk $\Delta$ in the complex plane,
defined by
\[
F_t(z)=1-(1-\exp(-t)+\exp(-t)\sqrt{1-z})^2, \quad z\in\Delta,\
t\ge0.
\]
One can check that $\Ss$ is generated by the following function:
\[
f(z)= -2\sqrt{1-z}(\sqrt{1-z}-1).
\]
It is easy to see that $f(1)=0$, but $f$ has no angular derivative
at the point $z=1$. At the same time, this point is not even a
fixed point of $\Ss$.

\vspace{2mm}

\noindent{\bf Example 5.} Now we consider the semigroup
$\Ss=\{F_t\}_{t\ge0}\subset\Hol(\Delta)$ defined by
\[
F_t(z)=\frac{(1+z)^{\alpha(t)}-(1-z)^{\alpha(t)}}{(1+z)^{\alpha(t)}+(1-z)^{\alpha(t)}}\,,
\]
where $\alpha(t)=e^{-2t}$. Differentiating at $t=0$, we find its
generator:
\[
f(z)= (1-z^2)\log{\frac{1+z}{1-z}}\,.
\]
Similarly as in the previous example, $f(\pm1)=0$, but $f$ has no
angular derivative at the points $\tau_1=1$ and $\tau_2=-1$.
However, in contrast with that example, these points are fixed
points of the semigroup:
\[
F_t(\pm1)=\pm1\quad\mbox{for all}\quad t\ge0.
\]
Moreover, it is possible to calculate the ``non-Euclidean
distance" $d_{\pm1}(F_t(z))$. In particular, for real $z=x$ we
have:
\[
d_1(F_t(x))=\frac{1-F_t(x)}{1+F_t(x)}=\left(\frac{1-x}{1+x}\right)^{\alpha(t)}=d_1(x)^{\alpha(t)}.
\]
Since $0<\alpha(t)<1$ when $t>0$, we conclude that for each fixed
$t>0$, the quotient $\displaystyle \frac{d_1(F_t(x))}{d_1(x)}$
tends to infinity as $z=x$ tends to $1$ radially, i. e.,
$d_1(F_t(z))$ does not admit an estimate of the form $A(t)d_1(z)$.

\bigskip

\noindent {\bf Acknowledgment.} The second author was partially
supported by the Fund for the Promotion of Research at the
Technion and by the Technion VPR Fund -- B. and G. Greenberg
Research Fund (Ottawa).

\vspace{\baselineskip}
\vspace{10mm}

\textit{Mark Elin}\newline Department of Mathematics,\newline ORT
Braude College,\newline 21982 Karmiel, Israel\newline
\textit{E-mail address: mark.elin@gmail.com}

\vspace{\baselineskip}
\vspace{10mm}

\textit{Simeon Reich}\newline Department of Mathematics,\newline
The Technion --- Israel Institute of Technology,\newline 32000
Haifa, Israel\newline \textit{E-mail address:
sreich@tx.technion.ac.il}

\vspace{\baselineskip}
\vspace{10mm}

\textit{David Shoikhet}\newline Department of Mathematics,\newline
ORT Braude College,\newline 21982 Karmiel, Israel\newline
\textit{E-mail address: davs27@netvision.net.il}

\end{document}